\documentclass[manmat,referee]{svjour}
\usepackage{graphics}
\usepackage{amsmath}
\usepackage{amssymb}
\usepackage{cite}
\newtheorem{notation}{Notation}

\numberwithin{equation}{section}
\begin{document}
\title{Pointwise Boundary Differentiability of Solutions of Elliptic Equations\thanks{Research supported by NSFC 11171266.}}
\author{Yongpan Huang \and Dongsheng Li \and Kai Zhang
}                     
%
%
\institute{School of Mathematics and Statistics, Xi'an Jiaotong University, Xi'an 710049, China.
 email: huangyongpan@gmail.com; lidsh@mail.xjtu.edu.cn; zkzkzk@stu.xjtu.edu.cn
\and Mathematics Subject Classification (2010): 35B65, 35J25, 35D40}
\date{Received: date / Revised version: date}
%
\maketitle
\begin{abstract}
In this paper, we give pointwise geometric conditions on the boundary which guarantee the differentiability of the solution at the boundary. Precisely, the geometric conditions are two parts: the proper blow up condition (see Definition \ref{d1.1}) and the exterior Dini hypersurface condition (see Definition \ref{d1.2}). If $\Omega$ satisfies this two conditions at $x_0\in\partial \Omega$, the solution is differentiable at $x_0$. Furthermore, counterexamples show that the conditions are optimal (see Remark \ref{re1.4} and the counterexample in Section 2).
\end{abstract}
\section{Introduction}
\label{}
We study the pointwise boundary differentiability of solutions of the following equations
\begin{equation}\label{e1.1}
 \left\{\begin{aligned}
 -a^{ij}(x)\frac{\partial^2u(x)}{\partial x_i \partial x_j}&=f(x)~~~~\mbox{in}~~\Omega; \\
 u&=0~~~~\mbox{on}~~\partial\Omega,
\end{aligned}
\right.
\end{equation}
where $\Omega\subset R^n~(n\geq 2)$ is a bounded domain; $(a^{ij}(x))_{n\times n}\in C(\bar{\Omega})$ is uniformly elliptic with $\lambda$ and $1/\lambda$ ($0<\lambda<1$) and $f\in C(\bar{\Omega})$. For convenience, solutions in this paper will always indicate viscosity solutions.

Li and Wang in \cite{LW} first obtained the boundary differentiability of the solutions of (\ref{e1.1}) for convex domains. Later, the same authors got the corresponding result in regard to the inhomogeneous boundary condition \cite{LW2}. Convexity is not the weakest geometric condition guaranteeing the boundary differentiability; Li and Zhang \cite{LZ} extended to a more general class of domains called $\gamma-$convex domains. As pointed out by the authors, $\gamma-$convexity is the weakest condition guaranteeing the boundary differentiability (see Remark 1.4\cite{LZ}).
However, \cite{LW}, \cite{LW2} and \cite{LZ} concerned the differentiability of the solution on the whole boundary, whereas we now focus on the pointwise situation.

%


For continuity of the solution up to the boundary, we only need an exterior restriction on the boundary such as the Wiener criterion for the boundary continuity (\cite{Wi}, see also p.330 \cite{Ke}), the exterior cone condition for the boundary H\"{o}lder continuity (p.205 \cite{GT}) and the exterior sphere condition for the boundary Lipschitz continuity (Lemma 1.2 \cite{Sa}) etc.. However, to obtain the continuity of the derivatives of the solution up to the boundary, an exterior restriction is not enough and we need the corresponding regularity on the boundary of the domain. For example, the boundary $C^{1,\alpha}$ regularity needs the $C^{1,\alpha}$ smoothness of the boundary (see \cite{Kr1}, \cite{Kr2} and \cite{MW}). The situation for the boundary differentiability is subtle. On the one hand, an exterior restriction is not enough and we need some other geometric properties of the boundary (see Definition \ref{d1.1} and the counterexample in Section 2). On the other hand, we don't need the $C^{1}$ (even differentiability) smoothness of the boundary. For instance, the boundary differentiability of the solution holds for convex domains whose boundary maybe not be differentiable at some points.

We give the following definition characterizing the geometric properties of the boundary.
\begin{definition}[\textbf{Proper blow up condition}]\label{d1.1}
Let $\Omega$ be a bounded domain and $x_0\in\partial\Omega$. We define the upper blow up set $C^{*}(x_0)$ and the lower blow up set $C_{*}(x_0)$ at $x_0$ as follows:
\begin{equation*}
    C^{*}(x_0):=conv\left(\{y:\forall~~ t>0,~~\exists~~ 0<t_1<t,~~\mbox{such that}~~x_0+t_1y\in\Omega\}\right),
\end{equation*}
and
\begin{equation*}
    C_{*}(x_0):=\{y:\exists~~ t>0,~~\forall~~ 0<t_1<t,~~x_0+t_1y\in\Omega\},~~~~~~~~~~~~~~~
\end{equation*}
where $conv$ denotes the convex hull. We say that $\Omega$ satisfies the proper blow up condition at $x_0\in\partial\Omega$ if
\begin{equation}
     C^{*}(x_0)\neq R^n_+ \mbox{ or } C_{*}(x_0)=R^n_+ \mbox{ (up to an isometry). }
\end{equation}

\end{definition}

\begin{remark}\label{re1.5}
We call $C$ a cone with vertex $x_0$ if $C$ is a convex set and $t(x-x_0)\in C$ for any $t>0$ and $x\in C$. It is called a proper cone if $C\neq R^n_+ $ (up to an isometry). By Definition \ref{d1.1}, $C^{\ast}$ is always a cone but $C_{\ast}$ may not be.
\end{remark}

Next, we introduce a definition for the exterior restriction on the boundary.
\begin{definition}[\textbf{Exterior Dini hypersurface condition}]\label{d1.2}
 We say that $\Omega$ satisfies an exterior Dini hypersurface condition at $x_0\in\partial\Omega$ if the following holds:
 there exists a unit vector $\eta\in R^n$, $R_0>0$ and $\gamma :R^+\longrightarrow R^+$ such that
 \begin{equation}\label{e1.2}
 \eta\cdot(x-x_0)\geq-\gamma(|x-x_0|),~\forall~x\in\bar\cap B(x_0,R_0),
 \end{equation}
 where $\gamma$ satisfies the Dini condition:
 \begin{equation}\label{e1.3}
 \int_{0}^{R_0}\frac{\gamma(r)}{r^2}dr<\infty.
 \end{equation}

\end{definition}

\begin{remark}\label{re1.6}
  (i) By \eqref{e1.3}, $\gamma(r)/r\rightarrow 0$ as $r\rightarrow 0$. Hence, $\gamma$ is differentiable at $0$ and $\gamma'(0)=0$.

  (ii) If $x_0=0$ and $\eta=e_n$, from \eqref{e1.2},
\begin{equation*}
  x_n\geq -\gamma (|x|).
\end{equation*}
Then $x_n/|x|\geq -\gamma(|x|)/|x| \rightarrow 0$ as $|x|\rightarrow 0$. Hence, there exist $r_1>0$ such that $x_n/|x'|>-1$ on $\bar{\Omega}\cap B_{r_1}$. Therefore,
\begin{equation}\label{e2.7}
  x_n\geq -\gamma(2|x'|)     \mbox{ on } \bar{\Omega}\cap B_{r_1},
\end{equation}
which means that the boundary representation function of $\Omega$ at $0$ lies above the hypersurface $-\gamma(2|x'|)$ locally. This is the geometrical explanation of the exterior Dini hypersurface condition.

(iii) If $\partial \Omega\in C^{1,\alpha}$ ($\alpha>0$), $\Omega$ satisfies the exterior Dini hypersurface condition at any boundary point.
\end{remark}

The following is our main result:
\begin{theorem}\label{th1.3}
Let $u$ be the solution of (\ref{e1.1}). If $\Omega$ satisfies the proper blow up condition and the exterior Dini hypersurface condition at $x_0\in\partial\Omega$, then $u$ is differentiable at $x_0$. That is, there exists a vector $\xi$ such that $u(x)=u(x_0)+\xi\cdot (x-x_0)+o(|x-x_0|),~\forall x\in\bar\Omega$.
\end{theorem}

\begin{remark}\label{re1.4}
(i) If $\Omega$ is convex, then for any $x_0\in\partial\Omega$, $C_{*}(x_0)=C^{*}(x_0)$ and $\Omega$ satisfies the exterior Dini hypersurface condition at $x_0$ with $\gamma\equiv 0$. Hence, the solution is differentiable at the whole boundary. This result was first obtained by Li and Wang \cite{LW}.

(ii) Furthermore, if $\Omega$ is $\gamma-$convex, then for any $x_0\in\partial\Omega$, $C_{*}(x_0)=C^{*}(x_0)$ and $\Omega$ satisfies the exterior Dini hypersurface condition at $x_0$ (see Definition 1.1 and Theorem 2.4 \cite{LZ}). Hence, the solution is differentiable at the whole boundary. Li and Zhang proved this in \cite{LZ}.

(iii) \cite{LW} and \cite{LZ} concerned the boundary differentiability on the whole boundary. For the pointwise boundary differentiability, Huang, Li and Wang \cite{HLW} proved that if $\partial\Omega$ is differentiable at $x_0$ and $\Omega$ satisfies the exterior Dini hypersurface condition there, then the solution of (\ref{e1.1}) is differentiable at $x_0$. This corresponds to the case $C^{*}(x_0)=C_{*}(x_0)=R^n_+$ in Theorem \ref{th1.3}.

(iv) \cite{LW}, \cite{LW2} and \cite{LZ} only concerned the case $C^{\ast}=C_{\ast}$. In fact, it is not necessary for the boundary differentiability.

(v) As pointed out in \cite{LZ} (Remark 1.4\cite{LZ}, see also Theorem 1.11\cite{Sa}), the exterior Dini hypersurface condition can not be weakened in Theorem \ref{th1.3}.

(vi) We call $\xi$ the gradient of $u$ at $x_0$ and write $Du(x_0)=\xi$.
\end{remark}

\begin{notation}
\begin{enumerate}
  \item For any $x\in R^{n}$, we may write $x=(x',x_{n})$, where $x'\in R^{n-1}$ and $x_{n}\in R$.
  \item $\{e_i\}^{n}_{i=1}$: the standard basis of $R^n$.
  \item $R^n_+:=\{x\in R^n|x_n>0\}.$
  \item $A^c:$ the complement of $A$ in $R^n,~\forall A\subset R^n$.
  \item $\bar A: $ the closure of $A,~\forall~A\subset R^n$.
  \item $\mbox{dist}(A,B):=\inf\{|x-y||x\in A, y\in B\}$.
  \item $a^{+}:=\max\{0,a\}$ \mbox{ and } $a^{-}:=-\min\{0,a\}$.
  \item $Q_r:=\{x:-r<x_i<r,~~\forall 1\leq i\leq n\}$ and $Q_r^+:=Q_r\cap R^n_+$.
\end{enumerate}
\end{notation}

\section{Differentiability at the boundary and a counterexample}
Li and Zhang \cite{LZ} proved the boundary differentiability for $\gamma-$convex domains which satisfy the exterior Dini hypersurface condition at any boundary point (see Definition 1.1 \cite{LZ}). For such a domain $\Omega$ and any boundary point $x_0$, the blow up set of $\Omega$ at $x_0$ is always a cone (i.e., $C^{\ast}(x_0)=C_{\ast}(x_0)$ is a cone, see Theorem 2.4 \cite{LZ}). If the cone is a half space, it is called a flat point. Otherwise, it is called a corner point (see Definition 2.5 \cite{LZ}). Li and Zhang proved the boundary differentiability separately with respect to these two kinds of boundary points. In the pointwise case, $C^{\ast}(x_0)\neq R^n_+$ and $C_{\ast}(x_0)= R^n_+$ are the counterparts of the corner point and the flat point respectively.

Let $\tilde{\gamma}: R^+\rightarrow R^+$ be a smooth function except at the origin and satisfy the Dini condition \eqref{e1.3} for some $R_0>0$. Let $\varphi $ be a function whose upper graph is a proper cone $C\subsetneq R^n_+$ with vertex the origin. Set
\begin{equation}\label{e2.6}
 \tilde{ \Omega}=:B_{R_0}\cap (\{x: x_n>\varphi(x') >0\}\cup \{x:x_n>-\tilde{\gamma}(|x'|) \mbox{ and } \varphi(x')=0 \}).
\end{equation}

The following two lemmas correspond to the two different kinds of boundary points.
\begin{lemma}\label{le2.1}
Let $u$ be the solution of (\ref{e1.1}) with $\Omega$ replaced by $\tilde{\Omega}$. Then $u$ is differentiable at the origin and $Du(0)=0$.
\end{lemma}
\noindent\textbf{Proof.} Clearly, the domain $\tilde{\Omega}$ satisfies the exterior Dini hypersurface condition with $\tilde{\gamma}$ at any boundary point. Hence, $\tilde{\Omega}$ is $\gamma-$convex. Furthermore, the blow up set of $\tilde{\Omega}$ at the origin is $C$ (the upper graph of $\varphi$) which is a proper cone. Hence, the origin is a corner point. By Theorem 3.3 \cite{LZ}, $u$ is differentiable at the origin and $Du(0)=0$.\qed


\begin{lemma}\label{le2.2}
Let $u$ be the solution of (\ref{e1.1}). Suppose that $\Omega$ satisfies the exterior Dini hypersurface condition at $x_0\in\partial \Omega$ and $\partial\Omega$ is differentiable at $x_0$. Then $u$ is differentiable at $x_0$.
\end{lemma}
\noindent\textbf{Proof.} Without loss of generality, we may assume that $x_0=0\in\partial\Omega$ and $\eta=e_n$ in Definition \ref{d1.2}. Since $\Omega$ satisfies the exterior Dini hypersurface condition and $\partial\Omega$ is differentiable at $0$, $C^{\ast}(0)=C_{\ast}(0)=R^n_+$. Hence, $0$ is a flat point. By Theorem 3.6 \cite{LZ}, $u$ is differentiable at $0$. In fact, the differentiability of $\partial \Omega$ at $0$ is the essence used in the proof of Theorem 3.6 \cite{LZ}.\qed

\begin{remark}\label{re1.7}
  We point out that Lemma \ref{le2.2} can be regarded as a special case of Theorem 2 \cite{HLW}.
\end{remark}


From above two lemmas, Theorem \ref{th1.3} follows easily.

\noindent\textbf{Proof of Theorem \ref{th1.3}.} Without loss of generality, we may assume that $x_0=0\in\partial\Omega$ and $\eta=e_n$ in Definition \ref{d1.2}. By the linearity of the equation, we assume that $u\geq 0$.

If $C^{*}(0)\neq R^n_+$, then by definition, $C^{*}(0)$ is a proper cone. On the other hand, for any $x_0\notin C^{\ast}(0)$, there exists $0<t_0<1$ such that $tx_0\in \Omega^c$ for all $0<t<t_0$. We choose another bigger proper cone $\tilde{C}$ with vertex the origin such that $C^{\ast}(0)\subset \tilde{C}$ and  $\Omega\cap B_{r_2}\cap R^n_+\subset \tilde{C}\cap B_{r_2}\cap R^n_+$ for some $r_2>0$. Next, we choose $\tilde{\gamma}:R^+\rightarrow R^+$ which is a smooth function except at the origin, satisfies the Dini condition \eqref{e1.3} and $\tilde{\gamma}(r)\geq \gamma(2r)$ for any $r>0$. Finally, construct the domain $\tilde{\Omega}$ as in \eqref{e2.6} with $R_0$ replaced by $r=\min(r_1,r_2)$ ($r_1$ is as in \eqref{e2.7}).  Then $\Omega\cap B_{r}\subset \tilde {\Omega}$. Let $\tilde {u}$
solve
\begin{equation*}\label{e2.8}
 \left\{\begin{aligned}
 -a^{ij}(x)\frac{\partial^2u(x)}{\partial x_i \partial x_j}&=f(x)~~~~\mbox{in}~~\tilde{\Omega}; \\
 u&=g~~~~\mbox{on}~~\partial\tilde{\Omega},
\end{aligned}
\right.
\end{equation*}
where $g\equiv 0$ on $\partial\tilde{\Omega}\cap \Omega^c$ and $g\equiv u$ on $\partial\tilde{\Omega}\cap \Omega$.

From Lemma \ref{le2.1}, $\tilde{u}$ is differentiable at $0$ and $D\tilde{u}(0)=0$. Since $\Omega\cap B_{r}\subset \tilde {\Omega}$, $0\leq u\leq \tilde{u} $ in $\Omega\cap B_{r}$. Hence, $u$ is differentiable at $0$ and $Du(0)=0$.

If $C_{*}(0)=R^n_+$, then by the definition, it is easy to know that $\partial \Omega$ is differentiable at $0$. From Lemma \ref{le2.2}, $u$ is differentiable at $0$. \qed\\

We have known that the exterior Dini hypersurface condition is necessary for boundary differentiability. Now we construct a counterexample to show that the proper blow up condition can not be dropped.

\begin{proposition}\label{pr2.3}
Let $\Omega\subset R^n$ ($n\geq 3$) be a bounded domain and $K\subset \Omega$ be compact whose $(n-2)$ dimension Hausdorff measure is zero, i.e., $H^{n-2}(K)=0$. If $\Delta u=0$ in $\Omega\backslash K$ and $u$ is bounded, then $K$ is removable, i.e., there exists an extension $\bar u$ of $u$ in $\Omega$ such that
\begin{equation*}
\left\{\begin{aligned}
\Delta \bar u&=0~~~~\mbox{in}~~\Omega;\\
\bar u&=u~~~~\mbox{in}~~\Omega\backslash K.
\end{aligned}\right.
\end{equation*}
\end{proposition}
\noindent\textbf{Proof.} See Theorem 1 (p.88) \cite{Ca} and Theorem 2 (p.151) \cite{EG}.\qed \\

From now on, we assume that the dimension $n\geq 4$. Denote $\tilde V_{x,\theta}:=\{y:\theta(y_n-x_n)\leq -|y'-x'|\}$ and $V_{x,\theta}:=|x|e_n+\tilde V_{x,\theta}$. From above proposition, we have the following lemma.
\begin{lemma}\label{le2.4}
Let $D_{\theta}=Q_1^+\backslash V_{1/4e_1,\theta}$ and $v^{\theta}$ solve
\begin{equation*}
\left\{\begin{aligned}
\Delta v^{\theta}&=0~~~~\mbox{in}~~D_{\theta};\\
v^{\theta}&=g^{\theta}~~~~\mbox{on}~~\partial D_{\theta},
\end{aligned}\right.
\end{equation*}
where $g^{\theta}\equiv 0$ on $\partial D_{\theta}\cap Q_1^+$ and $g^{\theta}\equiv x_n$ on $\partial D_{\theta}\backslash Q_1^+$. Then
\begin{equation}\label{e2.1}
    \|v^{\theta}-x_n\|_{L^{\infty}(D')}\rightarrow 0~~\mbox{as}~~\theta\rightarrow 0
\end{equation}
for any $D'\subset\subset \overline Q_1^+\backslash K_0$ where $K_0:=\{y:0\leq y_n\leq 1/4,y'=(1/4,0,...,0)\}$. Furthermore, we have
\begin{equation}\label{e2.2}
    \|v^{\theta}_n-1\|_{L^{\infty}(K_1)}\rightarrow 0~~\mbox{as}~~\theta\rightarrow 0,
\end{equation}
where $v^{\theta}_n:=\frac{\partial v^{\theta}}{\partial x_n}$ and $K_1:=\{y:y'=0,0\leq y_n\leq 1/2\}$.
\end{lemma}

\noindent\textbf{Proof.} Take the classical Schwarz reflection for $D^{\theta}$, $K_0$, $g^{\theta}$ and $v^{\theta}$, i.e.,
\begin{equation*}
    \tilde D_{\theta}=D_{\theta}\cup\{y:(y',-y_n)\in D_{\theta}\}\cup\{y:y_n=0,|y_i|<1,i\neq n,\mbox{ and }y\notin V_{1/4e_1,\theta}\},
\end{equation*}
$$\tilde K_0=K_0\cup\{y:(y',-y_n)\in K_0\},$$
\begin{equation*}
    \tilde g^{\theta}(x',x_n)=\left\{
    \begin{aligned}
     &g^{\theta}(x',x_n),~~~x_n\geq 0,\\
     &-g^{\theta}(x',-x_n),~~~x_n< 0
    \end{aligned}
    \right.
\end{equation*}
 and
\begin{equation*}
    \tilde v^{\theta}(x',x_n)=\left\{
    \begin{aligned}
     &v^{\theta}(x',x_n),~~~x_n\geq 0,\\
     &-v^{\theta}(x',-x_n),~~~x_n< 0.
    \end{aligned}
    \right.
\end{equation*}
Then $ \tilde v^{\theta}$ satisfies
\begin{equation*}
\left\{\begin{aligned}
\Delta \tilde v^{\theta}&=0~~~~\mbox{in}~~\tilde D_{\theta};\\
\tilde v^{\theta}&=\tilde g^{\theta}~~~~\mbox{on}~~\partial \tilde D_{\theta},
\end{aligned}\right.
\end{equation*}

By the maximum principle, $\|\tilde v^{\theta}\|_{L^{\infty}(\tilde D_{\theta})}\leq 1$. Then there exists a harmonic function $v$ in $Q_1\backslash \tilde K_0$ and a subsequence $\{v^{\theta_i}\}$ of $\{v^{\theta}\}$ such that $ \|v^{\theta_i}-v\|_{L^{\infty}(D')}\rightarrow 0$ as $i\rightarrow \infty$ for any $D'\subset\subset Q_1\backslash \tilde K_0$. Since $v^{\theta}$ is monotone (increasing or decreasing) for any $x\in D'$ as $\theta\rightarrow 0$, we have $ \|v^{\theta}-v\|_{L^{\infty}(D')}\rightarrow 0$. By Proposition \ref{pr2.3}, $\tilde K_0$ is removable. Therefore, $v\equiv x_n$ and (\ref{e2.1}) is proved.
By the interior estimates for the derivatives of the harmonic functions, we have
\begin{equation}\label{e2.3}
    \|v^{\theta}_n-1\|_{L^{\infty}(K_1)}\leq C\|v^{\theta}-x_n\|_{L^{\infty}(D'')}\rightarrow 0~~\mbox{as}~~\theta\rightarrow 0,
\end{equation}
where $D''\subset\subset Q_1\backslash \tilde K_0$ and $\mbox{dist}(K_1,\partial D'')>0$.
\qed\\

From (\ref{e2.2}), for $\theta_1$ small enough, we have
\begin{equation}\label{e2.4}
    v^{\theta_1}_n\mid_{K_1}\geq 1-1/4.
\end{equation}
Let $u^1=v^{\theta_1}$.


Similarly, we continue to construct $u^k$ for $k\geq 2$. Let $D_k=Q_1^+\backslash \bigcup_{m=1}^{k} V_{1/4^me_1,\theta_m}$ and $u^k$ solve
\begin{equation*}
\left\{\begin{aligned}
\Delta u^k&=0~~~~\mbox{in}~~D_k;\\
u^k&=g^k~~~~\mbox{on}~~\partial D_k,
\end{aligned}\right.
\end{equation*}
where $g^k\equiv 0$ on $\partial D_k\cap Q_1^+$ and $g^k\equiv x_n$ on $\partial D_k\backslash Q_1^+$. Then $u^k$ satisfies
\begin{equation}\label{e2.5}
    u^k_n\mid_{K_1}\geq 1-\sum_{m=1}^{k}1/4^m
\end{equation}
by choosing $\theta_m$ ($1\leq m\leq k$) properly.

Now we construct the counterexample.
\begin{theorem}[\textbf{A counterexample}]
Take
\begin{equation*}
    \Omega=Q_1^+\backslash\bigg(\bigcup_{m=1}^{\infty}V_{1/4^m,\theta_m}\bigg).
\end{equation*}
Let $g\equiv 0$ on $\partial\Omega\cap Q_1^+$ and $g\equiv x_n$ on $\partial\Omega\backslash Q_1^+$. Let $u$ solve the equation
\begin{equation*}
\left\{\begin{aligned}
\Delta u&=0~~~~\mbox{in}~~\Omega;\\
 u&=g~~~~\mbox{in}~~\partial\Omega.
\end{aligned}\right.
\end{equation*}
(Since $\Omega$ satisfies the exterior cone condition at any boundary point, the above equation is uniquely solvable). Then $u$ is not differentiable at $0$.
\end{theorem}
\begin{remark}
Obviously, $\Omega$ satisfies the exterior Dini hypersurface condition at $0$ with $\gamma\equiv 0$. But $C^{*}(0)=R^n_+$ and $C_{*}(0)\neq R^n_+$. We will show that $u$ is not differentiable at $0$. Hence, this is a counterexample to illustrate the necessariness of the proper blow up condition.
\end{remark}

\noindent\textbf{Proof.}  If $u$ is differentiable at $0$, it is easy to see that $\frac{\partial u}{\partial x_i}(0)=0$ for $1\leq i\leq n-1$. In addition, since $u(1/4^my_0)=0$ for any $m\geq 1$ where $y_0=(1,0,...,0,1)$, the directional derivative of $u$ along $y_0$ is zero. Note that $e_i$ ($1\leq i\leq n-1$) and $y_0$ are linear independent, thus $Du(0)=0$. In the following, we will obtain a contradiction with this.

By the maximum principle, $\|u^k-u\|_{L^{\infty}(\Omega)}\leq 1/4^k$. Thus, combining with (\ref{e2.5}), we have
\begin{equation*}
    u\mid K_1\geq u^k-1/4^k\geq (1-\sum_{m=1}^{k}1/4^m)x_n-1/4^k,~~\forall~~k\geq 1.
\end{equation*}
Therefore,
\begin{equation*}
    \frac{u(0,...,0,1/2^j)}{1/2^j}\geq (1-\sum_{m=1}^{j}1/4^m)-1/2^j\geq 1/6,~~\forall~~j\geq 1.
\end{equation*}
This contradicts with $Du(0)=0$.\qed\\

\bibliographystyle{amsplain}

\end{document}